\documentclass[journal]{IEEEtran}
\IEEEoverridecommandlockouts
\usepackage{cite}
\usepackage{amsmath,amssymb,amsfonts}
\AtBeginDocument{\pagenumbering{arabic}}

\usepackage{romannum}
\usepackage{graphicx}
\usepackage{mathtools}
\usepackage[percent]{overpic}  
\usepackage{svg}
\usepackage{xcolor}
\newtheorem{theorem}{Theorem}

\newtheorem{proposition}{Proposition}

\newtheorem{lemma}{Lemma}
\usepackage{array}
\usepackage{tabularx,booktabs}
\usepackage{booktabs,makecell,tabularx}

\usepackage{tikz}
\tikzset{
    block/.style = {rectangle, draw, text width=4cm, text centered, minimum height=1.5cm, node distance=1.5cm},
    smallblock/.style = {rectangle, draw, text width=3.2cm, text centered, minimum height=1cm, node distance=1.5cm},
    line/.style = {draw, -latex'},
}

\usepackage{pgfplots}
\usetikzlibrary{intersections}
\usepgfplotslibrary{fillbetween}
\usetikzlibrary{arrows, arrows.meta, positioning, backgrounds, automata, decorations, shapes, fit, calc}

\mathchardef\mhyphen="2D 

\ifCLASSOPTIONcompsoc
    \usepackage[caption=false, font=normalsize, labelfont=sf, textfont=sf]{subfig}
\else
\usepackage[caption=false, font=footnotesize]{subfig}
\fi

\usepackage{enumitem}

\def\bm#1{\boldsymbol{#1}}
\mathchardef\mhyphen="2D 

\newcommand{\cG}{{\cal G}}
\newcommand{\cW}{{\cal W}}

\newcommand{\cL}{{\cal L}}

\newcommand{\st}{\mbox{s.t.}}

\allowdisplaybreaks

\usepackage{hyperref}

\PassOptionsToPackage{hyphens}{url}\usepackage{hyperref}
\definecolor{lime}{HTML}{A6CE39}
\DeclareRobustCommand{\orcidicon}{%
	\begin{tikzpicture}
	\draw[lime, fill=lime] (0,0) 
	circle [radius=0.16] 
	node[white] {{\fontfamily{qag}\selectfont \tiny ID}};
	\draw[white, fill=white] (-0.0625,0.095) 
	circle [radius=0.007];
	\end{tikzpicture}
	\hspace{-2mm}
}

\foreach \x in {A, ..., Z}{%
	\expandafter\xdef\csname orcid\x\endcsname{\noexpand\href{https://orcid.org/\csname orcidauthor\x\endcsname}{\noexpand\orcidicon}}
}

\ifCLASSINFOpdf

\else

\fi

\pgfplotsset{compat=1.18}
\begin{document}

\title{FICA: Faster Inner Convex Approximation of Chance Constrained Grid Dispatch with Decision-Coupled Uncertainty}

\author{
        Yihong~Zhou\orcidA{},
        Hanbin~Yang\orcidB{},
        and Thomas~Morstyn\orcidC{},~\IEEEmembership{Senior Member,~IEEE} 
        \vspace{-5mm}
\thanks{Y. Zhou and  T. Morstyn are with the Department of Engineering Science, University of Oxford, U.K. (yihong.zhou@eng.ox.ac.uk, thomas.morstyn@eng.ox.ac.uk).
H. Yang is with the School of Data Science, The Chinese University of Hong Kong, Shenzhen, China (hanbinyang@link.cuhk.edu.cn).}
}

\markboth{IEEE Transactions on Power Systems}
{FICA: Faster Inner Convex Approximation of Chance Constrained Grid Dispatch with Decision-Coupled Uncertainty}

\maketitle

\begin{abstract}
This paper proposes a \textit{Faster Inner Convex Approximation} (FICA) method for solving power system dispatch problems with Wasserstein distributionally robust joint chance constraints (WJCC) and incorporating the modelling of the automatic generation control factors. The problem studied belongs to the computationally challenging class of WJCC with left-hand-side uncertainty (LHS-WJCC). By exploiting the special one-dimensional structure (even if only partially present) of the problem, the proposed FICA incorporates a set of strong valid inequalities to accelerate the solution process. We prove that FICA achieves the same optimality as the well-known conditional value-at-risk (CVaR) inner convex approximation method. Our numerical experiments demonstrate that the proposed FICA can yield $40\times$ computational speedup compared to CVaR, and can even reach up to $500\times$ speedup when the optimisation horizon exceeds 16 time steps. This speedup is achieved when only 50\% of constraints in a WJCC have the one-dimensional structure. The approximation quality is numerically verified to be the same as CVaR, and the quality gap is below 1\% when compared to the computationally demanding exact reformulation of the LHS-WJCC in most cases. We also discuss the applications of FICA in optimisation problems from other domains that (partially) exhibit the one-dimensional structure.

\looseness = -1
\end{abstract}



\IEEEpeerreviewmaketitle

\vspace{-4mm}
\section{Introduction}\label{sec:intro}

Power system dispatch represents a core class of decision-making problems in power system operation \cite{conejo2018power}. Its primary objective is to determine the power output of available generators in order to meet electricity demand at the lowest possible operational cost within technical constraints, and it is also the foundation of electricity market clearing \cite{chatzivasileiadis2018optimization}. However, the growing penetration of renewable generation is introducing new and significant challenges to such dispatch problems. In particular, the high uncertainty of renewable generation \cite{hong2016probabilistic} has made it increasingly important to adopt uncertainty-aware approaches to power system dispatch \cite{ordoudis2021energy, SFLA}.

An uncertainty-aware dispatch approach should account for both the reformulation of the objective function and the constraints. To handle uncertainty in the objective function, such as uncertain prices and the costs of corrective actions dependent on uncertain scenarios, existing studies typically reformulate the objective as a minimisation of the expected cost \cite{ming2017scenario}. This formulation aims to minimise the cost on average across a range of possible scenarios. However, power system operators place a strong emphasis on reliable system operation, and consequently, the management of uncertainty in the constraints becomes a critical focus, which is also the main focus of this paper. In the presence of uncertainty in the constraints, there are two main theoretical approaches. The first is to ensure constraint satisfaction under all possible (the worst-case) scenarios through robust optimisation (RO)---specifically, through the robust counterpart formulation \cite{wei2020tutorials}. The second is to guarantee constraint satisfaction with a high probability using chance constraints (CCs) \cite{shapiro2007tutorial}. In these two approaches, RO is often criticised for being overly conservative, as it hedges against the worst-case scenario \cite{kuhn2024distributionallyrobustoptimization, wei2020tutorials}. In contrast, CCs enable robust yet more economically efficient decision-making by explicitly limiting the probability of constraint violations. This provides a controlled level of conservativeness relative to RO, making CCs a less restrictive and more flexible alternative. In the context of power system dispatch, existing studies have extensively used CCs to ensure the high-probability satisfaction of both device-level and network-level constraints \cite{wu2016solution, yang2019analytical}.

Standard CCs rely on the assumption of exact knowledge of the uncertainty distribution. However, in practice, decision-makers typically have access only to a finite historical dataset. This dataset is often insufficient to accurately infer the true distribution of the random variables, thereby limiting the out-of-sample performance of standard CC methods \cite{mohajerin2018data}. This limitation has motivated the development of distributionally robust chance constraints (DRCCs) \cite{kuhn2024distributionallyrobustoptimization, mohajerin2018data}. DRCCs combine the principles of CCs and RO by ensuring the probability of constraint satisfaction under the worst-case distribution, which offers improved out-of-sample performance compared to standard CCs while still being less conservative than RO. In DRCCs, the worst-case distribution is defined with respect to a predefined ambiguity set, which can be constructed in various ways. Typical ambiguity sets can be defined by statistical moments (e.g., mean and variance) \cite{delage2010distributionally}, or by considering all distributions within a specified distance from a reference distribution. While moment-based ambiguity sets are computationally tractable, they can encompass an overly broad range of distributions, including those with significantly different and implausible shapes, resulting in overly conservative solutions \cite{gao2023distributionally, kuhn2024distributionallyrobustoptimization}. For distance-based ambiguity sets, typical distance measures include $\phi$-divergence (also known as $f$-divergence) and the Wasserstein distance \cite{kuhn2024distributionallyrobustoptimization}. The Wasserstein ambiguity set has gained widespread adoption due to its superior out-of-sample performance over other distance-based approaches \cite{mohajerin2018data, gao2023distributionally}, and has also been intensively investigated for power system dispatch problems \cite{chen2021fast, arrigo2022wasserstein}.

In a power system dispatch problem, there can be multiple constraints (e.g. for different generators, buses, and time steps) that involve uncertainty. In such case, a joint chance constraint (JCC) is preferred to guarantee simultaneous constraint satisfaction \cite{yzhou_agg, ding2022distributionally}. Combining with the Wasserstein ambiguity set leads to the Wasserstein distributionally robust joint chance constraint (WJCC). Unfortunately, WJCCs are challenging to solve, with the exact reformulation (where available) leading to mixed-integer programming (MIP) formulations with big-M coefficients, which are known to be computationally challenging and to scale poorly~\cite{ho2022distributionally, chen2024data, ho2023strong}. Recent studies~\cite{zhao2023viability} pointed out that these CC approaches are not practical for \emph{industry-scale} problems due to their computational intractability. This computational bottleneck has thus driven the development of approximate solution schemes for WJCC-based dispatch problems. For instance, the well-known CVaR approximation \cite{chen2023approximations} was employed in \cite{ordoudis2021energy}, while the iterative ALSO-X approximation algorithm \cite{jiang2022also} was used in \cite{wen2024multiple}.

Nonetheless, even with these approximation methods to solve the WJCC, the computational issue remains significant. The CVaR method introduces a large number of ancillary variables and constraints, which not only slows down computation but can also result in memory issues, as observed in our previous work \cite{SFLA}. Iterative approximation methods---such as ALSO-X \cite{jiang2022also} and the subsequent ALSO-X\# \cite{jiang2024also}, although demonstrating better optimality than CVaR, have been shown to require several times or even several tens of times more computational time compared to CVaR in different case studies \cite{jiang2022also, jiang2024also}. As pointed out by \cite{bukhsh2024balancing}, real-world power system dispatch---such as the Great Britain Balancing Mechanism---involves many complex dispatch requirements, rendering the problem computationally intensive and mixed-integer even in a deterministic setting. Computational efficiency is therefore a critical factor for an operational framework to be used by system operators \cite{bukhsh2024balancing}. With the increasing integration of controllable devices into the grid (such as storage), computationally efficient solutions to WJCC are urgently needed to enable the real-world application of WJCC-based dispatch.

Being aware of this critical computational issue, our previous work \cite{SFLA} proposed a strengthened and faster linear approximation (SFLA) to WJCC, demonstrating up to a $100\times$ speedup compared to the CVaR approach. However, SFLA is specifically designed for WJCCs with right-hand-side uncertainty (RHS-WJCC) \cite{ho2022distributionally}, which limits its applicability to many power system dispatch problems that must account for re-dispatch actions, such as those involving automatic generation control (AGC) factors \cite{ning2021deep} (the AGC factors are also referred to as participation factors in existing studies~\cite{geng2019data}). These actions necessitate a left-hand-side formulation (LHS-WJCC), because the generator output power will depend on the product of their assigned AGC factor and the system imbalance \cite{ho2023strong, ning2021deep}. In fact, the LHS-WJCC is known to be computationally challenging, and its general form \emph{cannot} even be exactly reformulated as a MIP problem \cite{ho2023strong}. However, we observe that the LHS-WJCC structure in power system dispatch, arising from AGC factors, (partially) exhibits a special one-dimensional structure that can be leveraged to extend our previous success with the SFLA for RHS-WJCC. Motivated by this insight, we propose a \emph{Faster Inner Convex Approximation} (FICA) for this type of special LHS-WJCC, achieving significant computational speedup. In summary, this paper makes the following contributions:
\begin{enumerate}
    \item We propose the \emph{Faster Inner Convex Approximation} (FICA), a method that exploits the one-dimensional structure of LHS-WJCCs in power system dispatch problems involving AGC factors. By harnessing this structure, FICA achieves up to $500\times$ speedup compared to the well-known CVaR method.
    \item We theoretically prove that despite the significant computational speedup, FICA maintains equivalent optimality to the CVaR approach.
    \item Numerical experiments and comparisons with benchmarks across a wide range of parameter settings demonstrate the effectiveness of the proposed FICA and its robust computational advantages.
\end{enumerate}

The remainder of this paper is organised as follows. Section \ref{sec:WJCCED} presents the fundamental mathematical formulation of the WJCC-equipped power system dispatch problems. Section \ref{sec:FICA} explores the special structure of these problems and introduces the proposed FICA, with its theoretical equivalence to CVaR analysed in Section \ref{sec:eq-cvar}. Section \ref{sec:casestudies} provides numerical case studies, Section \ref{sec:final_remark} discusses the broader impact of FICA, and Section \ref{sec:conclusion} concludes this paper.

\section{Joint Chance Constrained Economic Dispatch}\label{sec:WJCCED}

Let $\cL$ denote the set of transmission lines, $\cG$ collect the indices of dispatchable generators, and $\cW$ contain the wind farms. For each wind farm $w \in \cW$ at time step $t \in [T]\coloneqq\{1,\ldots,T\}$, its wind generation is uncertain and can be modelled as the sum of a deterministic forecast $\omega_{t,w}$ and a random forecasting error $e_{t,w}$. Following the data-driven manner, we assume that we have $N$ historical data samples for the random vector $\bm e \coloneqq [\bm e_t]_{t\in[T]}$ with $\bm e_t \coloneqq [e_{t, w}]_{w\in\cW}$ (this vector definition helps capture correlations). We use $[N] \coloneqq \{1, \cdots, N\}$ to index the samples of a random vector throughout the paper, such that we have the dataset as $\{ \bm e_i \}_{i\in[N]}$ (without the tilde overhead). We also assume that the network demand $\bm d_t$ is deterministic (considering its uncertainty would not affect our subsequent derivation). Following \cite{geng2019data} and \cite{ning2021deep}, we consider the following power system dispatch problem:
\begin{subequations} \label{prob:ed}
    \begin{align}
        \min_{\bm p, \bm\alpha} \ \ & \mathcal{C}(\bm p, \bm \alpha) \label{obj:dispatch}\\
        \st \ \ & \bm p,\ \bm\alpha \in \mathcal{X}, \label{cons:ed-general}\\
        & \bm 1^\top (\bm p_t - \bm d_t) = 0, &&\hspace{-71pt} \forall t\in[T], \label{cons:ed-balance} \\
        &  \bm 1^\top \bm \alpha_t = 1, &&\hspace{-71pt} \forall t\in[T], \label{cons:ed-agc} \\
        &  -\bm 1 \leq \bm \alpha_t \leq \bm 1, &&\hspace{-71pt} \forall t\in[T], \label{cons:ed-agc-range} \\
        & \underline{\bm p}_t \leq \bm p_t \leq \overline{\bm p}_t, &&\hspace{-71pt} \forall t\in[T], \label{cons:ed-power-sch} \\
         \inf\limits_{\mathbb{P} \in \mathcal{F}(\theta)} \!\! & \mathbb{P}\! \left\{ 
        \begin{array}{@{}l@{}}
            \underline{\bm p}_t \leq \bm p_t - \! \sum\limits_{w\in\cW} e_{t,w} \bm \alpha_t \leq \overline{\bm p}_t, \hspace{14pt} \forall t\in[T], \\
            \underline{\bm f} \leq 
            {\bm S^G} \tilde{\bm p}_t +\\
            \ \ \quad {\bm S^W} (\bm \omega_t + \bm e_t)-{\bm S^D} \bm d_t \leq \overline{\bm f}, \forall t\in[T]
        \end{array}
    \right\}\!\! \geq 1 - \epsilon, \label{cons:ed-jcc}
    \end{align}
\end{subequations}
where $\tilde{\bm p}_t \coloneqq \bm p_t - \sum_{w\in\cW} e_{t,w} \bm \alpha_t$ represents the actual generation power after AGC adjustments.
The objective \eqref{obj:dispatch} is to minimise the total dispatch cost $\mathcal{C}(\bm p, \bm \alpha)$, which is a function of the generation dispatch $\bm p \coloneqq (\bm p_t)_{t \in [T]}$ with $\bm{p}_t \coloneqq (p_{t,g})_{g \in \cG}$, and the AGC factors $\bm\alpha \coloneqq (\bm\alpha_t)_{t \in [T]}$ with $\bm{\alpha}_t \coloneqq (\alpha_{t,g})_{g \in \cG}$. Constraint \eqref{cons:ed-general} allows the decision variables $\bm p$ and $\bm\alpha$ to be subject to additional arbitrary constraints contained in $\mathcal{X}$, which can include ramp rate constraints or even non-convex unit commitment constraints. Constraint \eqref{cons:ed-balance} ensures the power balance of the dispatch schedule. Constraint \eqref{cons:ed-agc} ensures that the corrective actions determined by the AGC factors $\bm \alpha$ can fully compensate for the system imbalance caused by forecasting errors. Constraint \eqref{cons:ed-agc-range} allows negative AGC values following \cite{geng2019data}, which can be helpful to alleviate network congestion. Constraint \eqref{cons:ed-power-sch} ensures that the scheduled dispatch $\bm p$ remains within the power limits of all generators.

The large WJCC \eqref{cons:ed-jcc} ensures that constraints with uncertain parameters are jointly satisfied for all generators, transmission lines, and time steps, with a probability not less than $1-\epsilon$ under the worst-case probability distribution in a Wasserstein ambiguity set $\mathcal{F}(\theta)$. For a random vector $\bm\xi \in \mathbb{R}^K$, $\mathcal{F}(\theta)$ collects all distributions within a Wasserstein $\theta$-radius ball centred at an empirical distribution $\mathbb{P}_N$ (from historical data): 
\begin{align*}
    \mathcal{F}(\theta) \coloneqq \left\{ \mathbb{P} \,\middle|\, \inf_{\Pi \in \mathcal{P}(\mathbb{P}_N, \mathbb{P})} \mathbb{E}_{(\boldsymbol{\xi}, \boldsymbol{\xi}') \sim \Pi} \left[ \| \boldsymbol{\xi} - \boldsymbol{\xi}' \| \right] \leq \theta \right\},
\end{align*}
where $\mathcal{P}(\mathbb{P}_N, \mathbb{P})$ is a set of distributions on $\mathbb{R}^K \times \mathbb{R}^K$ with marginal distributions $\mathbb{P}_N$ on $\mathbb{R}^K$ and $\mathbb{P}$ on $\mathbb{R}^K$, following ~\cite{ho2022distributionally, chen2024data, ho2023strong}. Let $[P]$ denote the index set of all individual constraints involved in the large WJCC~\eqref{cons:ed-jcc}. These constraints can be categorised into two groups: the subset $[P]^*$ corresponds to the first line in~\eqref{cons:ed-jcc}, which ensures that the actual generation power remains within the power limits of the generators, while the complement set $[P] \setminus [P]^*$ corresponds to the second line in~\eqref{cons:ed-jcc}, which ensures that the line flow limits are respected under uncertain wind generation $(\bm\omega_t + \bm e_t)$ and the adjusted generation $\tilde{\bm p}_t$. In this second group, a direct-current power flow (DCPF) model is applied, with $\bm S^G$, $\bm S^W$, and $\bm S^D$ representing the corresponding power transfer distribution factor (PTDF) matrices \cite{chatzivasileiadis2018optimization}.

\section{FICA: Faster Inner Convex Approximation}\label{sec:FICA}
\subsection{Challenges of the General LHS-WJCC}\label{subsec:challenge}
    
    The WJCC \eqref{cons:ed-jcc} in the power system dispatch problem is an LHS-WJCC, because the AGC adjusted generation $\tilde{\bm p}_t \coloneqq \bm p_t - \sum_{w\in\cW} e_{t,w} \bm \alpha_t$ involves a multiplication between decision variables $\bm \alpha$ and random variables $\bm e$, and this term appears in both groups of constraints in \eqref{cons:ed-jcc}. Following the convention in the existing literature \cite{chen2024data, ho2023strong}, an LHS-WJCC consists of the following standard form of constraints:
    \begin{equation}\label{eq:gel-lhs-wdrjcc}
        (\bm b_p - \bm A_p^\top \bm x)^\top \bm \xi + d_p - \bm a_p^\top \bm x \geq 0, \quad p\in [P],
    \end{equation}
    where $\bm b_p$, $\bm A_p$, $d_p$, and $\bm a_p$ denote the coefficients associated with each of the $P$ constraints in an LHS-WJCC. Here, $\bm x$ is the decision variable vector, and $\bm \xi$ is the random vector. However, solving this class of LHS-WJCC remains highly challenging. As discussed in~\cite{ho2023strong}, an exact MIP reformulation is only available under a homogeneity assumption, i.e., when $\bm b_p = \bm b$ and $\bm A_p = \bm A$ for all $p \in [P]$. While it is theoretically possible to enforce homogeneity by introducing artificial degeneracy into the distribution (e.g., expanding $\bm A_p$ with redundant rows), doing so often introduces extra conservativeness. Furthermore, the reformulated problem still leads to a large-scale MIP with big-M constraints, which are known to be computationally intensive and to scale poorly with the number of constraints or decision variables.

\subsection{Analysis of the Special Structure}\label{sec:analysis_1dim}
    Fortunately, the LHS-WJCC~\eqref{cons:ed-jcc} has a special structure that can be exploited for computational advantages. Rather than involving a multiplication between a decision vector and a random vector as in the general LHS-WJCC \eqref{eq:gel-lhs-wdrjcc}, the LHS term in the $[P]^*$ group of the WJCC \eqref{cons:ed-jcc} (i.e., the generation power limits) consists of 
    \emph{only one} multiplication between a one-dimensional decision variable $\alpha_{t,g}$ and a one-dimensional random variable $\sum_{w\in\cW} e_{t,w}$. Specifically, for each constraint $p \in [P]^*$, the LHS term takes the following form:
    \begin{equation}\label{eq:1d-structure-compact}
        (\bm b_p - \bm A_p^\top \bm x)^\top \bm \xi = \alpha_{t,g} \sum\limits_{w\in\cW} e_{t,w} \mbox{ or } - \alpha_{t,g} \sum\limits_{w\in\cW} e_{t,w},
    \end{equation}
    which corresponds to the upper or lower generation limits $\overline{p}_{t,g}$ and $\underline{p}_{t,g}$, respectively. This one-dimensional structure leads to the following theoretical property:
    \begin{proposition}\label{prop:1dim}
        Suppose $\sum\limits_{w\in\cW} {e}_{1,t,w} \leq \cdots \leq \sum\limits_{w\in\cW} {e}_{N,t,w}$, then the $j$-th smallest element of $\{\alpha_{t,g} \sum\limits_{w\in\cW} {e}_{i,t,w} \}_{i\in[N]}$ is equal to $\min \{ \alpha_{t,g} \sum\limits_{w\in\cW} {e}_{j,t,w},\  \alpha_{t,g} \sum\limits_{w\in\cW} {e}_{N-j+1,t,w} \}$.
    \end{proposition}
    The rationale is that the relative order of the sequence elements is reversed when $\alpha_{t,g}$ is negative, and preserved when $\alpha_{t,g}$ is positive. 
    Note that only the $[P]^*$ set of constraints in the WJCC~\eqref{cons:ed-jcc} has the aforementioned one-dimensional structure. In contrast, the constraints in $[P]\setminus [P]^*$ \emph{adds an extra random variable term} $\bm S^W \bm e_t$ to $\alpha_{t,g} \sum\nolimits_{w\in\cW} {e}_{i,t,w}$. This destroys the one-dimensional structure, and the beneficial structural property cannot be leveraged. However, as will be demonstrated, even partial exploitation of the one-dimensional structure in $[P]^*$ can significantly improve the computational efficiency of solving the LHS-WJCC~\eqref{cons:ed-jcc}.


    
\subsection{The Proposed FICA Formulation}

    The proposed FICA extends the success of our previous SFLA in RHS-WJCC \cite{SFLA}. An RHS-WJCC consists of the following set of constraints:
    \begin{equation*}
        \bm b_p^\top \bm \xi + d_p - \bm a_p^\top \bm x \geq 0, \  p\in [P].
    \end{equation*}
    We can observe that the only difference between the RHS-WJCC and the LHS-WJCC is the coefficient of $\bm \xi$. Based on this observation, it is clear that for each specific value of the decision vector $\bm x$, a general LHS-WJCC reduces to a RHS-WJCC since the term $(\bm b_p - \bm A_p^\top \bm x)$ becomes a parameter. This motivates us to directly apply the SFLA formulation \cite{SFLA} to the LHS-WJCC as follows:
    \begin{subequations}\label{eq:SFLA}
	    \begin{align}
	        & s \geq 0, \bm r \geq \bm 0, \label{eq:SFLA_s_r} \\
	        & \epsilon N s - \sum\limits_{i\in [N]}r_i \geq \theta N, \label{eq:SFLA_epsilonNs} \\
	        & \kappa_i \Big(\dfrac{(\bm b_p - \bm A_p^\top \bm x)^\top \bm \xi_i + d_p - \bm a_p^\top \bm x}{\| \bm b_p - \bm A_p^\top \bm x \|_*} \Big) \geq s - r_i,\nonumber \span\\
            & & \forall i \in [N]_p,\ p\in [P],  \label{eq:SFLA_main} \\
	        & \dfrac{q_p(\bm x) + d_p - \bm a_p^\top \bm x}{\| \bm b_p - \bm A_p^\top \bm x \|_*} \geq s, & \forall p \in [P], \label{eq:SFLA_spup}
        \end{align}
    \end{subequations}
    where $s$ and $\bm r$ are ancillary decision variables. The term $\bm \xi_i$ with $i\in[N]$ represents the $i$-th data sample for the random vector $\bm \xi$ (namely $\bm e$ in the dispatch problem \eqref{prob:ed}). The dual norm $\|\cdot\|_*$ corresponds to the norm used in defining the Wasserstein distance. Vector $\bm \kappa \coloneqq [\kappa_i]_{i \in [N]}$ represents tunable hyperparameters, and a uniform setting $\bm \kappa = \bm 1$ is recommended in~\cite{SFLA}. 
    Both $q_p(\bm x)$ and $[N]_p$ play a key role in enhancing computational efficiency by tightening the space of ancillary variables and removing redundant constraints, which are defined as:
    \begin{align*}
    	& q_p(\bm x) \coloneqq \text{($k{+}1$)-th smallest value of } \{ (\bm b_p - \bm A_p^\top \bm x)^\top \bm \xi_i \}_{i\in[N]}, \\
    	& [N]_p \coloneqq \{i \in [N] \mid (\bm b_p - \bm A_p^\top \bm x)^\top \bm \xi_i < q_p(\bm x) \}.
    \end{align*}
    Although a general LHS-WJCC can be expressed using SFLA as in \eqref{eq:SFLA}, such a formulation is not usable due to the difficulty of obtaining a closed-form expression for $q_p(\bm x)$ and the constraints in $[N]_p$. However, an exception arises in the one-dimensional special structure of the LHS-WJCC \eqref{cons:ed-jcc} in the power system dispatch problem. This special structure provides us with an analytical expression for $q_p(\bm x)$ illustrated in Proposition \ref{prop:1dim}. Specifically, for any $p \in [P]^*$, we have:
    {\small\begin{align}\label{eq:q_1dim}
        q_p(\bm x) = & \min\{ (\bm b_p - \bm A_p^\top \bm x)^\top \bm \xi_{k+1}, (\bm b_p - \bm A_p^\top \bm x)^\top \bm \xi_{N-k} \},
    \end{align}}where, without loss of generality, we assume that the indices $1, \ldots, N$ are ordered such that
    \begin{equation}\label{eq:order}
        (\bm b_p - \bm A_p^\top \bm 1)^\top \bm \xi_1 \leq \cdots \leq (\bm b_p - \bm A_p^\top \bm 1)^\top \bm \xi_N.
    \end{equation}
    These elements are derived by setting $\bm x = \bm 1$, which is the same as setting $\alpha_{t,g}=1$ in Proposition~\ref{prop:1dim}. 
    When constructing constraints, the $\min(\cdot)$ operator in~\eqref{eq:q_1dim} can be replaced by two inequality constraints, each corresponding to one of the two elements inside $\min(\cdot)$. This is valid since $q_p(\bm x)$ appears on the left-hand side of a ``$\geq$" constraint. Similarly, the set $[N]_p$ can be expressed as:
    \begin{align}\label{eqn:reductionset}
        [N]_p = \{ 1,\cdots, k\} \cup \{N-k+1, \cdots, N \}.
    \end{align}
    It includes the indices of the first $k$ and the last $k$ elements 
    of the ordered sequence in Eq.~\eqref{eq:order}
    , which correspond to the two terms inside the $\min(\cdot)$ operator.

    Since only the $[P]^*$ group of the WJCC \eqref{cons:ed-jcc} enjoys this one-dimensional structure, the SFLA can only be partially adapted, resulting in the following migrated version:
	\begin{subequations}\label{eq:SFLA-mig}
	    \begin{align}
	        & s \geq 0, \bm r \geq \bm 0,  \label{eq:SFLA-mig_s_r} \\
	        & \epsilon N s - \sum\limits_{i\in [N]}r_i \geq \theta N,  \label{eq:SFLA-mig_epsilonNs} \\
	        & \kappa_i \Big(\dfrac{(\bm b_p - \bm A_p^\top \bm x)^\top \bm \xi_i + d_p - \bm a_p^\top \bm x}{\| \bm b_p - \bm A_p^\top \bm x \|_*} \Big) \geq s - r_i,  \nonumber \span\\ 
            & &\forall i \in [N]_p,\ p\in [P]^*,  \label{eq:SFLA-mig_main} \\
	        & \dfrac{q_p(\bm x) + d_p - \bm a_p^\top \bm x}{\| \bm b_p - \bm A_p^\top \bm x \|_*} \geq s, & \forall p \in [P]^*,  \label{eq:SFLA-mig_spup}\\
	        & \kappa_i \Big(\dfrac{(\bm b_p - \bm A_p^\top \bm x)^\top \bm \xi_i + d_p - \bm a_p^\top \bm x}{\| \bm b_p - \bm A_p^\top \bm x \|_*} \Big) \geq s - r_i,  \nonumber \span \\ 
            & & \forall i \in [N],\ p\in [P] \setminus [P]^*.  \label{eq:SFLA-mig_main2}
	    \end{align}
	\end{subequations}
    The form in \eqref{eq:SFLA-mig} is non-convex because the decision variable $\bm x$ appears in the denominator of constraints \eqref{eq:SFLA-mig_main}--\eqref{eq:SFLA-mig_spup}. To mitigate this issue, we propose a further inner approximation by replacing the denominator $\| \bm b_p - \bm A_p^\top \bm x \|_*$ with $\max_{p\in[P]}\| \bm b_p - \bm A_p^\top \bm x \|_*$. After this replacement, the denominator becomes uniform across all constraints, which motivates us to multiply both sides by this common denominator. This leads to the proposed \textit{Faster Inner Convex Approximation} (FICA) expressed in the following general form:
    \begin{subequations}\label{eq:FICA}
	    \begin{align}
	        & s \geq 0, \bm r \geq \bm 0, \label{eq:FICA_s_r} \span \\
	        & \epsilon N s - \sum\limits_{i\in [N]}r_i \geq \theta N \max_{p\in[P]}\| \bm b_p - \bm A_p^\top \bm x \|_*, \label{eq:FICA_epsilonNs} \span\\
	        & \kappa_i \Big((\bm b_p - \bm A_p^\top \bm x)^\top \bm \xi_i + d_p - \bm a_p^\top \bm x \Big) \geq s - r_i, \notag \span\\
            & \qquad\qquad\qquad\qquad\qquad\qquad \forall i \in [N]_p,\ p\in [P]^*,  \label{eq:FICA_main} \span \\
	        & \kappa_i \Big((\bm b_p - \bm A_p^\top \bm x)^\top \bm \xi_i + d_p - \bm a_p^\top \bm x \Big) \geq s - r_i, \notag \span \\
            & \qquad\qquad\qquad\qquad\qquad\quad \forall i \in [N],\ p\in [P] \setminus [P]^*,\label{eq:FICA_main2} \span \\
	        & (\bm b_p - \bm A_p^\top \bm x)^\top \bm \xi_{k+1} + d_p - \bm a_p^\top \bm x \geq s, & \forall p \in [P]^*, \label{eq:FICA_s_str_0} \span \\
            & (\bm b_p - \bm A_p^\top \bm x)^\top \bm \xi_{N-k} + d_p - \bm a_p^\top \bm x \geq s, & \forall p \in [P]^*, \label{eq:FICA_s_str} \span
	    \end{align}
	\end{subequations}
    where constraints \eqref{eq:FICA_s_str_0} and \eqref{eq:FICA_s_str} correspond to \eqref{eq:SFLA-mig_spup} that expands the $\min(\cdot)$ operator in \eqref{eq:q_1dim}.
	Appendix~\ref{appen:fica-ed} presents the FICA formulation using the specific nomenclature of the power system dispatch problem~\eqref{prob:ed}.

\section{Equivalence to CVaR}\label{sec:eq-cvar}

An alternative and widely applied inner approximation is the CVaR approach \cite{mohajerin2018data, chen2023approximations}, which can be formulated as the following set of convex constraints:
\begin{subequations}
    \begin{align*}
        & \bm \alpha \geq \bm 0, \beta \in \mathbb{R}, \tau \in \mathbb{R}, \\
        & \tau + \frac{1}{\epsilon} \Big(\theta \beta + \frac{1}{N} \sum_{i\in [N]} \alpha_i\Big) \leq 0,  \span\\
        & \alpha_i \geq w_p \Big(\bm a_p^\top \bm x - (\bm b_p - \bm A_p^\top \bm x)^\top \bm \xi_i - d_p\Big) - \tau,  \span\\
        & & \forall i\in[N], p\in [P], \\
        & \beta \geq w_p \| \bm b_p - \bm A_p^\top \bm x \|_*,  & \forall p \in [P], 
    \end{align*}
\end{subequations}
where $\bm w \coloneqq [w_p]_{p\in[P]} \in \{\bm w \in (0,1)^P \mid \sum_{p\in[P]}w_p = 1\}$ are tunable hyperparameters that affect the performance of the CVaR approximation by prioritising specific constraints within a WJCC. The tunable weights $\bm w$ are proposed by \cite{chen2023approximations}, with standard CVaR having $w_p = 1/P$ for each $p\in[P]$. In the following, we will show that CVaR with $\bm w = \bm 1/P$ has a \emph{projected} feasible region for $\bm x$ that is equivalent to that of the proposed FICA method, given the recommended setting $\bm \kappa=\bm 1$. Our analysis focuses on $\bm x$-space, as the remaining ancillary variables do not appear in the objective function.

Note that for any optimal solution to the CVaR formulation, we can always set $\beta = \max_{p\in[P]} w_p \| \bm b_p - \bm A_p^\top \bm x \|_*$ without affecting optimality (see the proofs of Proposition 11 in \cite{chen2023approximations}). This allows us to eliminate the variable $\beta$ from the CVaR expression as follows:
\begin{subequations}
    \begin{align}
        & \bm \alpha \geq \bm 0, \tau \in \mathbb{R},  \\
        & \tau + \frac{1}{\epsilon} \Big(\theta \max_{p\in[P]} w_p \| \bm b_p - \bm A_p^\top \bm x \|_* + \frac{1}{N} \sum_{i\in [N]} \alpha_i\Big) \leq 0,   \\
        & \alpha_i \geq w_p \Big(\bm a_p^\top \bm x - (\bm b_p - \bm A_p^\top \bm x)^\top \bm \xi_i - d_p\Big) - \tau,  \notag \\
        & && \hspace{-97pt}  \forall i\in[N], p\in [P]. 
    \end{align}
\end{subequations}
Let $s = -\tau$ and $\bm r = \bm\alpha$, we reach the following equivalent form for CVaR with the standard $\bm w = \bm 1/P$:
\begin{subequations}\label{eq:cvar-lhs-eq}
    \begin{align}
        & \bm r \geq \bm 0, s \geq 0,  \label{eq:cvar-lhs-eq-s-r}\\
        & \epsilon N s - \sum_{i\in [N]} r_i \geq \theta N \max_{p\in[P]} \| \bm b_p - \bm A_p^\top \bm x \|_* ,   \label{eq:cvar-lhs-eq-epsNs}\\
        & (\bm b_p - \bm A_p^\top \bm x)^\top \bm \xi_i - \bm a_p^\top \bm x + d_p \geq s - r_i,   \notag \\ & && \hspace{-97pt} \forall i\in[N], p\in [P], \label{eq:cvar-lhs-eq-main}
    \end{align}
\end{subequations}
Denote the $\bm x$-\emph{feasible} region of CVaR by:
\begin{align*}
    \mathcal{X}_C \coloneqq \{\bm x \mid \exists s, \bm r \ \text{satisfying \eqref{eq:cvar-lhs-eq-s-r}--\eqref{eq:cvar-lhs-eq-main} }\},
\end{align*}
which is the projection of of the feasible region in~\eqref{eq:cvar-lhs-eq} onto the $\bm x$-space.
We then have:
\begin{theorem}\label{theorem:cvar_subeq_fica_exist}
    For any $\bm x \in \mathcal{X}_C$, there exists ($\bm r$, $s$) such that constraints~\eqref{eq:FICA_s_r}--\eqref{eq:FICA_s_str} are satisfied, under $\bm \kappa = \bm 1$.
\end{theorem}
\begin{theorem}\label{theorem:fica_subeq_cvar_all}
    For any $(\bm x, \bm r, s)$ that satisfies the FICA constraints in \eqref{eq:FICA} under $\bm \kappa = \bm 1$, constraints~\eqref{eq:cvar-lhs-eq-s-r}--\eqref{eq:cvar-lhs-eq-main} that define $\mathcal{X}_C$ are satisfied.
\end{theorem}
The proofs are provided in Appendix \ref{appen:proof}. The two theorems above establish the equivalence of the $\bm x$-feasible regions between the proposed FICA and the standard CVaR approximation. Furthermore, the stronger result in Theorem \ref{theorem:fica_subeq_cvar_all}---which holds for $\bm x, \bm r, s$ rather than just $\bm x$ as in Theorem \ref{theorem:cvar_subeq_fica_exist}---implies that the proposed FICA reduces the feasible space for the ancillary variables $\bm r$ and $s$. This reduction is the key to the computational speedup achieved by the proposed FICA.

\section{Numerical Case Studies}\label{sec:casestudies}

\subsection{Experiment Settings}\label{sec:settings}
\begin{figure}[tb]
    \centering
    \vskip -3mm
    \subfloat[ ]{\includegraphics[width=.9\linewidth]{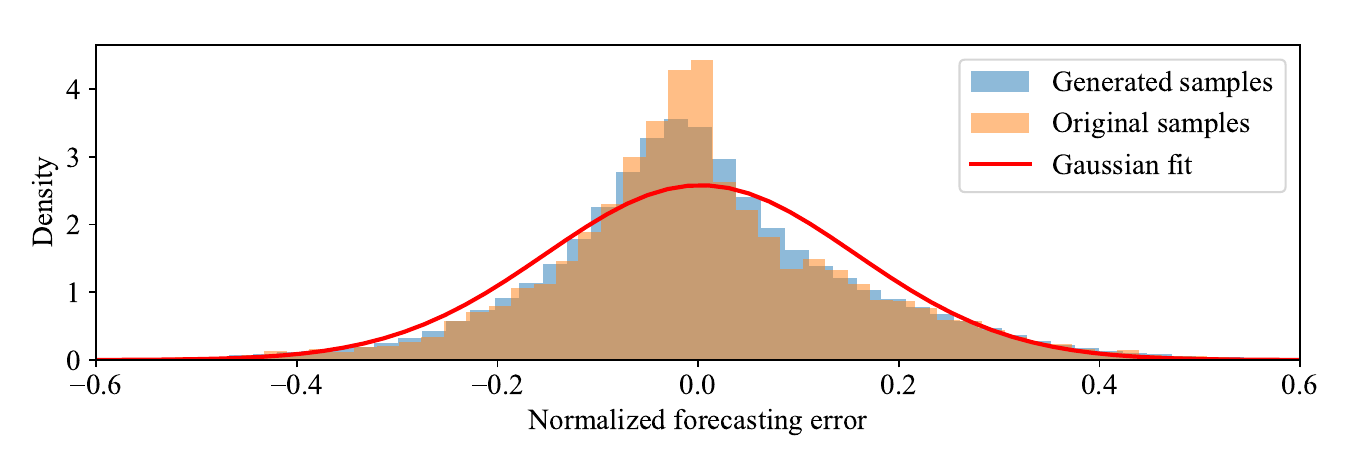}\label{fig:marginal}}\\
    \subfloat[ ]{\includegraphics[width=1\linewidth]{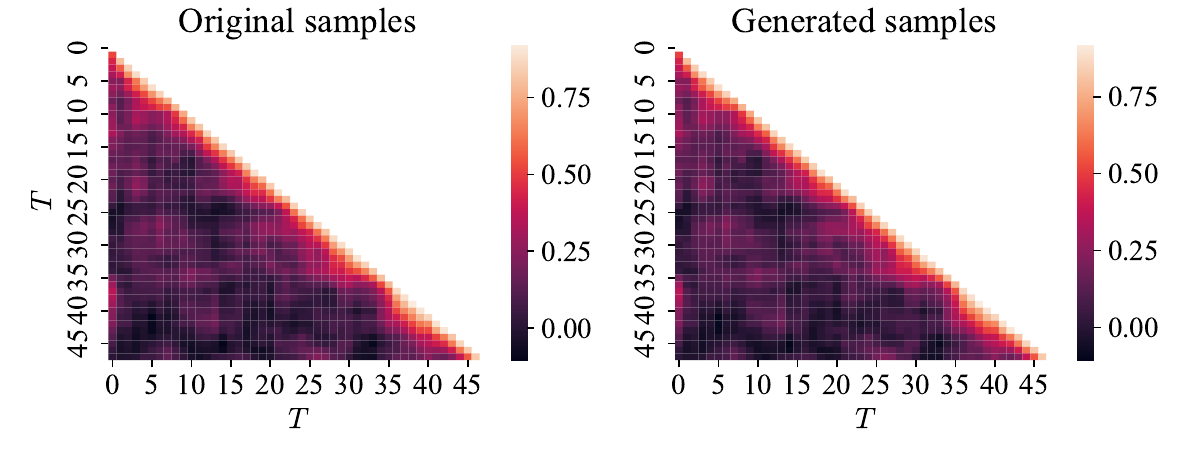}\label{fig:corr}}
    \vskip -2mm
    \caption{Comparison of original samples and generated samples for wind forecasting error. Original samples are those 156 forecasting error samples from Global Energy Forecasting Competition 2012. Generated samples refer to those generated by our fitted KDE plus copula.  (a) The marginal distribution of the wind forecasting error normalised by the wind generator capacity. We plot a Gaussian fit to highlight the non-Gaussian shape of our samples. (b) The temporal Pearson correlation coefficients.}
    \vskip -2mm
    \label{fig:wind_dist}
\end{figure}
\begin{figure}[tb]
    \centering
    \vskip -3mm
    \includegraphics[width=1\linewidth]{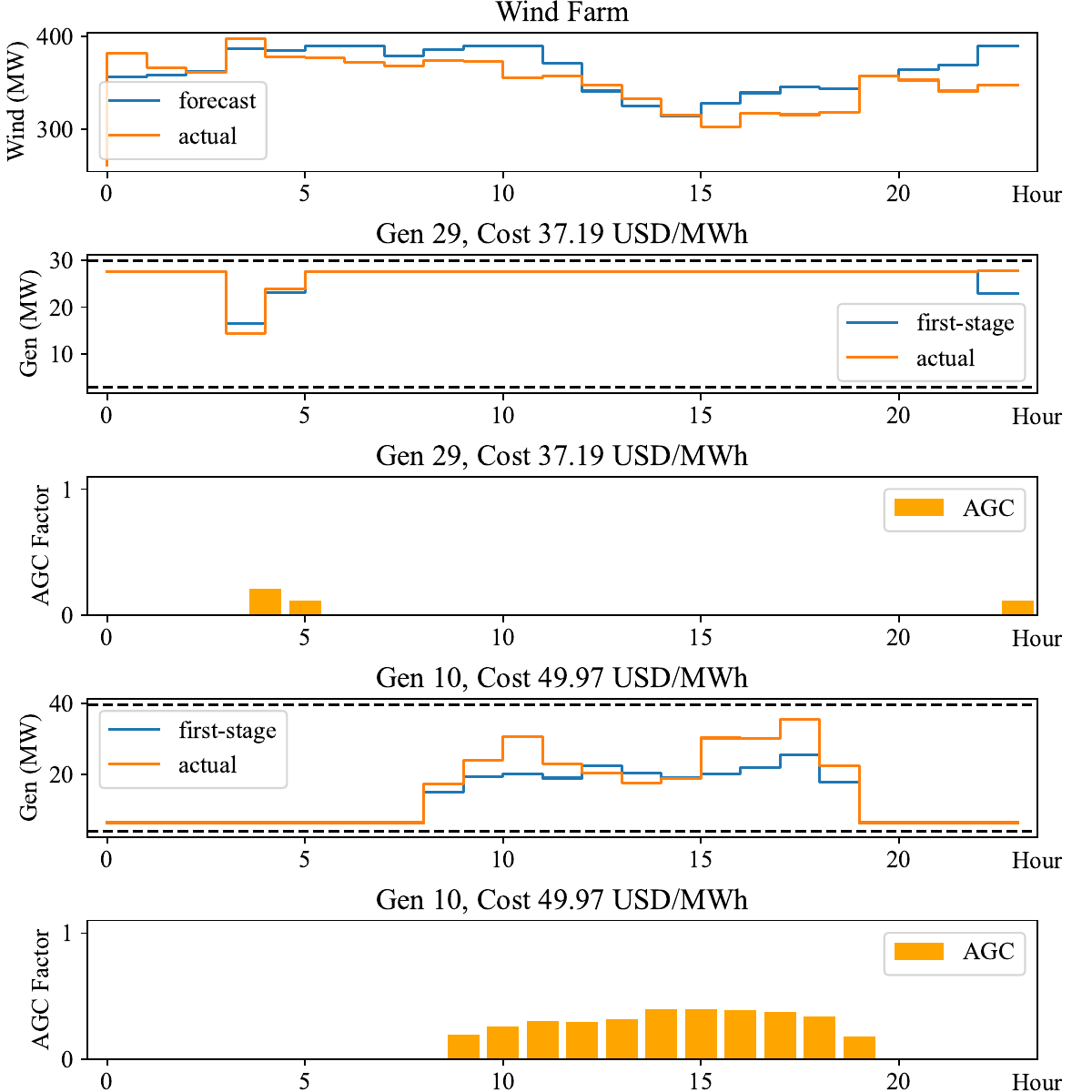}
    \vskip -3mm
    \caption{Day-ahead wind forecast and generator power profiles versus the actual wind power and the adjusted generator power alongside the AGC factors. Dashed black lines represent the generator power limits.}
    \vskip -2mm
    \label{fig:dispatch}
\end{figure}
\begin{figure*}[tb]
    \centering
    \vskip -3mm
    \includegraphics[width=1\linewidth]{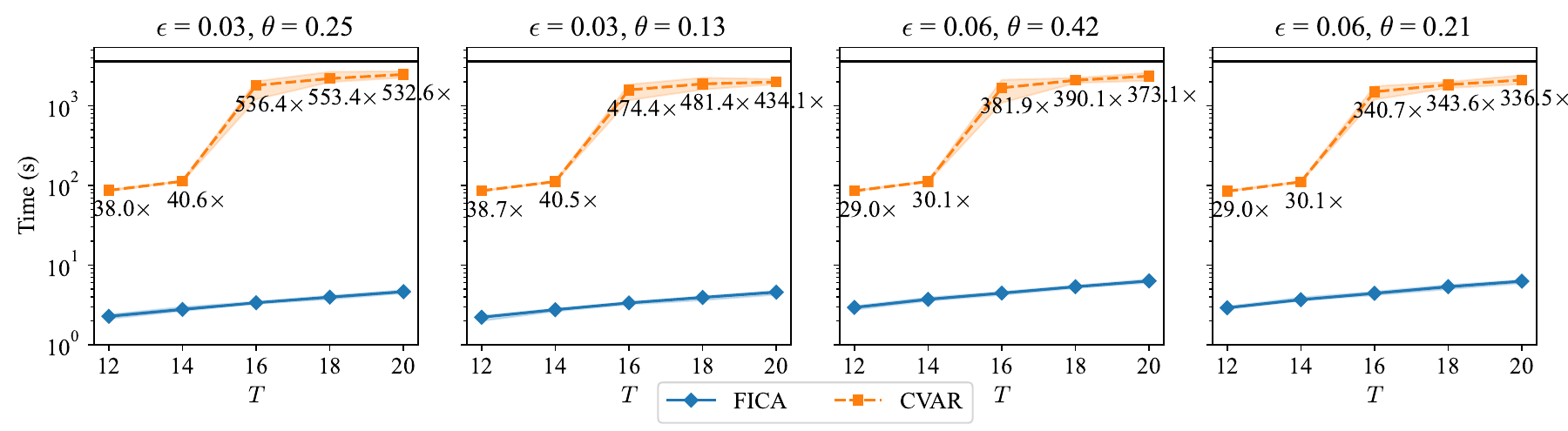}
    \vskip -3mm
    \caption{Comparison of the computing time between the proposed FICA and CVaR under different settings of risk levels $\epsilon$ and Wasserstein radii $\theta$. Dots represent the mean values over the ten random runs, while shaded areas represent the 99\% percentile intervals. Numbers ending with a $\times$ indicate the speedup rate (CVaR time / FICA time).}
    \vskip -3mm
    \label{fig:time_large}
\end{figure*}

This section presents numerical case studies to demonstrate the computational speedup and the approximation quality of the proposed FICA. The dispatch problem is implemented on the standard IEEE 24-bus test case with 38 branches \cite{ieee24bus}, using its default load and generation settings. We consider a quadratic cost function as the dispatch objective~\eqref{obj:dispatch}, which is independent of the AGC factor $\bm\alpha$ and the random variables $\bm e$, namely,
\begin{equation}
    \mathcal{C}(\bm p, \bm\alpha) \equiv \mathcal{C}(\bm p) = \sum\nolimits_{t\in[T], g\in\cG} (\bar{\bar{c}}_g p_{t,g}^2 + \bar{c}_g p_{t,g} + c_g ).
\end{equation}
The generator cost coefficients ($\bar{\bar{c}}_g$, $\bar{c}_g$, $c_g$) and other technical parameters follow \cite{xu2017application}.
This alleviates the need to formulate an expectation minimisation, as our key focus is on the LHS-WJCC \eqref{cons:ed-jcc}. Expectation minimisation under the Wasserstein metric can be reformulated as a set of convex constraints \cite{ordoudis2021energy} that can be incorporated into $\mathcal{X}$ in problem \eqref{prob:ed}, so the proposed FICA can be applied as well. We add $10$ wind farms to the network with a total capacity of $730$\,MW, representing $60\%$ of the total network load and a high level of renewable penetration.

Following \cite{yang2019analytical}, we consider uncertain wind power generation, modelled as a deterministic forecast plus a random forecasting error. Wind data, including deterministic forecasts and error samples, are obtained from the publicly available Bronze-medal forecasts (by Mudit Gaur) for the second wind farm in the Global Energy Forecasting Competition 2012 \cite{gefc2012}. The competition task was to forecast wind power generation for the next $48$ hours based on real-world wind farm data. We use the out-of-sample forecasts from the Bronze-medal results and compute the difference to the released ground-truth data, obtaining $156$ samples of wind power forecasting error, each spanning $48$ hours. These samples enable us to calculate the temporal correlation and the marginal distributions. To enrich the sample size, we use a kernel density estimate (KDE) to fit the marginals and apply the Gaussian copula to generate temporally correlated samples that respect the original marginal distribution. Fig.~\ref{fig:wind_dist} compares the original $156$ samples with those generated by the KDE plus copula, demonstrating a good fit that preserves both the non-Gaussian marginals and the original temporal correlation.

Both of the numerical cases are implemented in \text{Python 3.9.19} and solved by \text{Gurobi 11.0.3} on a server with dual AMD EPYC 7532 CPUs and 256 GB RAM. The solver is configured as \text{FeasibilityTol=$10^{-9}$}, \text{OptimalityTol=$10^{-9}$}, \text{IntFeasTol=$10^{-9}$}, \text{Threads=4}, and a one-hour wall-clock time limit \text{TimeLimit=3600s}; other solver parameters are set to default values. We set $\bm \kappa=\bm 1$ for the proposed FICA and $\bm w=\bm 1/P$ for CVaR following \cite{ordoudis2021energy}. For all norm calculations in the WJCC formulations, we adopt the $L_1$ norm, and its dual is the $L_\infty$ norm. To obtain robust conclusions, we conduct $10$ independent random runs for each parameter setting. In each random run, we 1) sample a distinct set of generator cost parameters based on the range suggested by~\cite{xu2017application}, which impacts the objective function; 2) sample a different set of historical samples $\{ \bm \xi_i \}_{i \in [N]}$ for the random wind forecasting error, which affects the feasible region; 3) select a random starting simulation hour, uniformly drawn from the $24$ hours of a day, which affects the demand and generation profiles and thereby the feasible region; 4) set a different seed for the optimisation solver, affecting the solver behaviors such as tie-breaking rules. Finally, to evaluate the satisfaction probability (reliability) of the WJCC, we test the solution on $5{,}000$ out-of-sample scenarios. Codes are publicly available at our GitHub repository: \url{https://github.com/PSALOxford/FICA}. 

\subsection{Dispatch Results}\label{sec:dispatch}

We begin with a sanity check by showcasing the dispatch results obtained after solving the WJCC-equipped power system dispatch \eqref{prob:ed}. As shown in Fig.~\ref{fig:dispatch}, generator $29$ is scheduled to operate close to its maximum output for most of the time, with a slight margin to hedge against the distributional ambiguity. This is due to its lower generation cost (the linear component of the quadratic cost). The higher-cost generator $10$ is kept close to its minimum generation except when it is assigned AGC factors, in which periods the day-ahead scheduled generation is kept at a mid-level to mitigate the risk of being unable to adjust in real time. Finally, these generators adjust their power output in response to wind shortages and surpluses as expected. These results verify the correctness of our model formulation.

\subsection{Computational Speedup}\label{sec:case_time}

\begin{figure}[tb]
    \centering
    \vskip -2mm
    \includegraphics[width=1\linewidth]{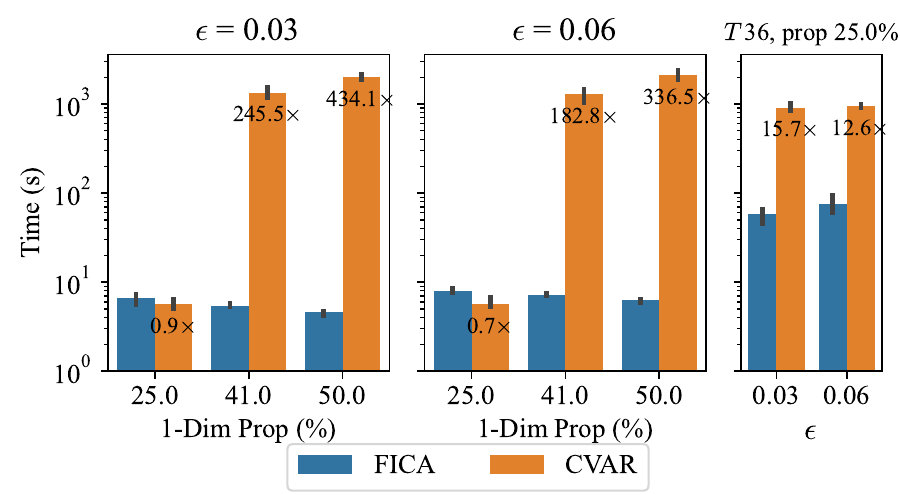}
    \vskip -4mm
    \caption{Comparison of computing time between the proposed FICA and the well-known CVaR. \textit{1-Dim Prop} represents the percentage proportion of the one-dimensional-structure constraints $[P]^*$ over the whole constraints in WJCC \eqref{cons:ed-jcc}. The rightmost subplot shows the computing time when the problem scales up for $\textit{1-Dim Prop}=25\%$. Numbers ending with a $\times$ indicate the speedup rate (CVaR time / FICA time). The short black vertical lines on top of the bars represent the $99\%$ percentile intervals. $\epsilon$ is the risk level.}
    \vskip -3mm
    \label{fig:time_onedim_prop}
\end{figure}

Fig.~\ref{fig:time_large} shows the computing time required to solve the proposed FICA and the well-known CVaR method in \eqref{eq:cvar-lhs-eq}. The number of generators is set to $38$, which is the same as the number of network branches. For each risk level $\epsilon$, we have tested two values of the Wasserstein radius $\theta$, one leading to slightly higher out-of-sample reliability and the other to slightly lower, to mimic the tuning process of $\theta$. As can be seen, when the optimisation horizon $T$ is less than or equal to $14$, the proposed FICA achieves a $30\times$ to $40\times$ speedup. A noteworthy observation is the sharp increase in computing time for CVaR when $T \geq 16$, in which the proposed FICA achieves a substantial speedup of $300\times$ to $500\times$. Note that this section does not compare the computing time against the exact MIP formulation of LHS-WJCC given by~\cite{ho2023strong}, as it fails to find a feasible solution within the one-hour time limit here.

As discussed in Section \ref{sec:FICA}, the successful application of FICA relies on the one-dimensional structure of the generation power limit constraints within the large-scale LHS-WJCC \eqref{cons:ed-jcc}. To understand how this structural feature impacts performance, we analyze the sensitivity of FICA’s computational speedup to the proportion of one-dimensional constraints among all WJCC constraints, referred to as \textit{1-Dim Prop}. Fig.~\ref{fig:time_onedim_prop} displays these results by changing the number of generators. The two subplots on the left indicate that as \textit{1-Dim Prop} increases, FICA indeed leads to greater computational speedup. When \textit{1-Dim Prop} is low ($25\%$), CVaR is even slightly faster, with a $10$--$30\%$ reduction in computing time. However, it should be noted that in these \textit{1-Dim Prop} $=25\%$ cases, the computing time for both methods is below $10$ seconds, making this slight inferiority of FICA minor. Furthermore, even under the \textit{1-Dim Prop} $=25\%$ setting, when the problem scales up (e.g., $T=36$ in the rightmost subplot of Fig.~\ref{fig:time_onedim_prop}), the proposed FICA regains its computational advantage and achieves more than a $10\times$ speedup. Combining these results with the substantial speedup for large $T$ observed in Fig.~\ref{fig:time_large}, we can conclude that the computational advantage of the proposed FICA becomes significant for large-scale problems. This is particularly the case in power system dispatch, where the large $T$, the substantial number of dispatchable generators, and the large-scale nature of power networks all hinder the application of WJCC. Here, the proposed FICA may bring a significant computational advantage in these industry-scale problems, thereby unclocking the value of WJCC. Finally, we have also observed significant memory consumption when using CVaR. In contrast, our FICA is more memory-saving, because its $[N]_p$ sets eliminate a huge portion of constraints in CVaR. As can be seen in Eq.~\eqref{eqn:reductionset}, the $[N]_p$ set for each $p\in [P]^*$, which consists of two parts, only have $2\lfloor \epsilon N \rfloor$ constraints, as opposed to the $N$ constraints for each $p\in [P]$ for CVaR.

\subsection{Validation of FICA’s Optimality}\label{sec:case_opt}
\begin{figure}[tb]
    \centering
    \vskip -2mm
    \includegraphics[width=1\linewidth]{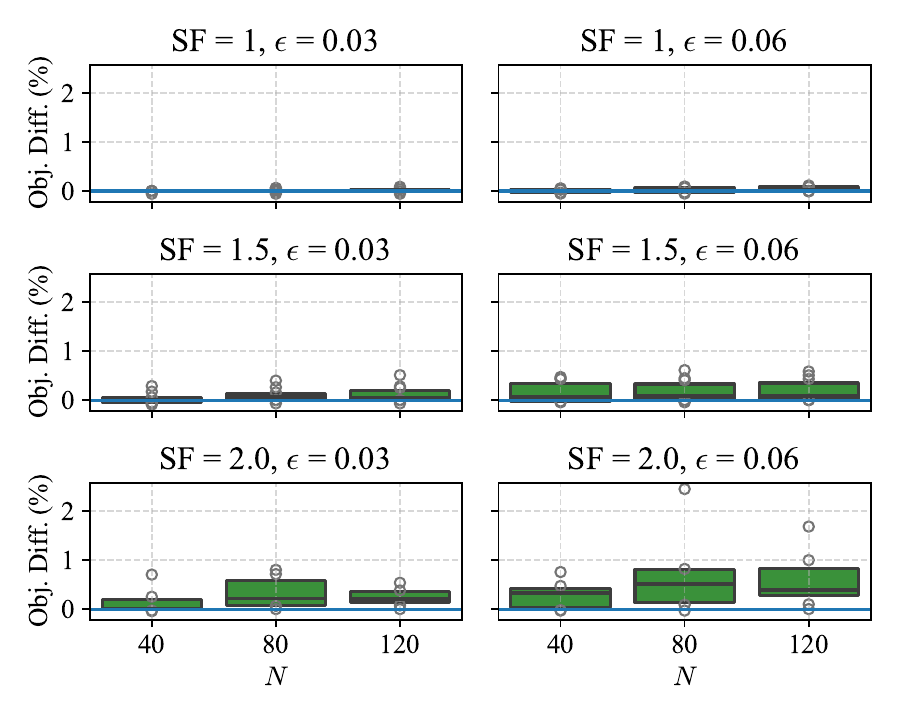}
    \vskip -6 mm
    \caption{Percentage difference of the optimal objective values obtained by the proposed FICA and the exact solution method. Each box collects the percentage difference for the ten random runs. A positive value indicates that the proposed FICA has worse optimality. The blue lines highlight the places with zero differences. SF stands for the load scaling factor that scales up the network load, which then increases network congestion and creates harder cases. $\epsilon$ is the risk level, and $N$ is the number of scenarios used for WJCC.}
    \label{fig:optimality}
\end{figure}

Fig.~\ref{fig:optimality} compares the optimality of the proposed FICA with the exact reformulation method (MIP) for LHS-WJCC, where artificial degeneracy is added for the exact reformulation following \cite{ho2023strong}, otherwise no MIP reformulation exists as discussed in Section \ref{subsec:challenge}. We applied the formulation (11) in \cite{ho2023strong} instead of their (25) as it requires excessive extra computing time to solve $NP$ subproblems only for reducing the big-M values in (11). The optimisation horizon is set to $1$ in this section, given the excessive computing difficulty of the exact solution method. The comparison with CVaR is not included, as we have verified that FICA achieves the same optimality (with differences less than $10^{-7}$\%) as CVaR for all test cases. We have observed that the tightness of a JCC affects the approximation quality, as the feasible region to be approximated becomes more complex. Therefore, we perform the comparison across different risk levels $\epsilon$, different numbers of scenarios $N$, and different load scaling factors (SF)---which scale up the network load, increase network congestion, and thus increase the tightness of the JCC. In Fig.~\ref{fig:optimality}, a positive value (\textit{Obj. Diff}) indicates that the exact method can further reduce the objective value by a certain percentage of the objective value solved by FICA. As can be seen, this \textit{Obj. Diff} is below $1\%$ for most cases, with only one random run exceeding $2\%$. Overall, we can conclude that the proposed FICA achieves great approximation quality in power system dispatch problems.

\section{Discussions}\label{sec:final_remark}
\vspace{-1mm}

Some final remarks are provided to highlight the features of the proposed FICA. First, the proposed FICA is applicable when there are several WJCCs in a single problem, which enables more targeted risk management. The theoretical properties remain valid, as one can incorporate other JCCs into the general feasible region $\mathcal{X}$ \eqref{cons:ed-general} when proving Theorems \ref{theorem:cvar_subeq_fica_exist} and \ref{theorem:fica_subeq_cvar_all}. Moreover, the proposed FICA imposes no restriction on the convexity of the original problem, unlike the scenario-based JCC approximation which requires convexity \cite{geng2019data}. 

Our proposed FICA enables the partial exploitation of the special one-dimensional structure of constraints by introducing the flexible $[P]^*$ set without any restrictions in deriving the theoretical results. It has demonstrated strong computational speedup through such partial exploitation. This capability is important, as it is very likely that a single JCC contains constraints with different structures, where only part of them have specific computational properties. Furthermore, since FICA generalises our previous SFLA, this partial-exploitation property also holds for the combination of FICA and SFLA. In other words, one can combine the strengthening constraints from FICA and SFLA (such as \eqref{eq:FICA_s_str}) to accelerate the solution of a WJCC that contains both RHS constraints and one-dimensional LHS constraints, without compromising optimality compared to CVaR.

The one-dimensional structure is not exclusive to power system dispatch problems. In the context of bidding problems in power system reserve markets, this structure arises in modelling the activated reserve capacity, which is represented as the product of the reserve capacity (i.e., bids) and the uncertain reserve activation proportion. Moreover, the one-dimensional structure is commonly encountered in other domains. In finance, it appears in risk-sharing problems, where the objective is to determine the optimal proportion of risk (e.g., financial losses) shared among peers. In logistics, a similar structure emerges when modelling the optimal assignment proportion of demand to different supply centres, as a way to hedge against deviations in customer demand.

In addition, the current dispatch formulation \eqref{prob:ed} is closely related to multi-stage stochastic programming (MSSP) problems. An MSSP problem aims to determine the optimal first-stage decision that minimises the expected cost while accounting for the ability to make subsequent recourse actions as uncertainty is gradually revealed. Our studied dispatch formulation \eqref{prob:ed} mirrors this paradigm by identifying a one-dimensional linear decision rule $\bm \alpha$ to conservatively approximate the optimal recourse actions. Consequently, the proposed FICA can also be employed as an inner-approximation to a WJCC-equipped MSSP. Such formulations are crucial to ensure that prior-stage decisions can have ``safe'' future-stage recourse actions with a high probability. Although this paradigm has been established for decades \cite{prekopa2013stochastic_programming}, it has seen limited application due to its computational complexity. Ref.~\cite{ordoudis2021energy} studied WJCC-equipped MSSP in a two-stage setting with a linear decision rule, solved using CVaR. Given the superior performance of the proposed FICA over CVaR, we believe that FICA unlocks more opportunities to realise the practical value of WJCC-equipped MSSP.

Recent advances for solving LHS-WJCC include several iterative solution schemes, such as ALSO-X~\cite{jiang2022also}, ALSO-X\#~\cite{jiang2024also}, and the ``Terminator''~\cite{jiang2025terminator}. We did not compare the proposed FICA with these approaches because either 1) their case studies have demonstrated longer computing times than CVaR~\cite{jiang2022also, jiang2024also}, and/or 2) their methods use CVaR as an initial solution and then further improve approximation optimality \cite{jiang2024also, jiang2025terminator}. The significant computational speedup of the proposed FICA over CVaR thus nullifies the need for a direct comparison of computing time with these methods. The key advantage of these iterative algorithms \cite{jiang2022also, jiang2024also, jiang2025terminator} lies in their ability to improve the optimality of CVaR to near-optimal levels, very close to the exact WJCC solution. However, in our power system dispatch problems, the optimality gap of the proposed FICA is already less than $1\%$ compared to the exact method in most cases. Even in scenarios where FICA or CVaR do not provide high-quality solutions, the proposed FICA can be seamlessly integrated into these iterative algorithms due to its proven equivalence to CVaR.

\section{Conclusion}\label{sec:conclusion}

This paper proposes a \textit{Faster Inner Convex Approximation} (FICA) for solving WJCC-equipped power system dispatch problems that include AGC factor modelling. The problem studied falls within the challenging class of LHS-WJCCs. However, by exploiting the special one-dimensional structure, the proposed FICA extends our previous success in the simpler RHS-WJCC setting. We proved that FICA always achieves the same optimality as the well-known convex inner-approximation, CVaR. Our numerical experiments demonstrate that the proposed FICA can achieve up to $500\times$ computational speedup compared to CVaR and with the same optimality. Finally, compared to the computationally intensive exact reformulation of LHS-WJCC, the approximation quality of FICA is only compromised by less than $1\%$ for most cases. 

{
\setlength{\abovedisplayskip}{2pt} 
\setlength{\belowdisplayskip}{2pt} 
\appendices

\section{FICA for Power System Dispatch} \label{appen:fica-ed}

We re-express the key terms in FICA \eqref{eq:FICA} under the nomenclature for the power system dispatch problem \eqref{prob:ed}. The following content divides the constraints in the WJCC \eqref{cons:ed-jcc} into those for $\overline{\bm p}_t$, $\underline{\bm p}_t$, $\overline{\bm f}_t$, and $\underline{\bm f}_t$, respectively. Note that the provided FICA form is for the recommended $\bm \kappa = \bm 1$. And we assume 
\[\sum\limits_{w\in\cW} {e}_{1,t,w} \leq \cdots \leq \sum\limits_{w\in\cW} {e}_{N,t,w},\]
which is equivalent to~\eqref{eq:order}.

The following content details the terms of FICA \eqref{eq:FICA} that need to be transformed into the nomenclature for the power system dispatch problem \eqref{prob:ed}: 
\subsubsection*{For the upper generation limit  $\overline{\bm p}_t$ constraints}
\begin{align*}
    & \bm b_p - \bm A_p^\top \bm x = \alpha_{t,g}\bm 1 \in \mathbb{R}^{|\cW|}, \\
    & (\bm b_p \! - \! \bm A_p^\top \bm x)^\top \bm\xi_i \! + \! d_p \! - \! \bm a_p^\top \bm x = \overline{p}_{t,g} + \alpha_{t,g}(\sum\limits_{w\in \cW}e_{i,t,w}) - p_{t,g}, \\
    & q_p(\bm x) + d_p - \bm a_p^\top \bm x = \\
    & \overline{p}_{t,g} + \min\{ \alpha_{t,g}(\sum\limits_{w\in \cW}e_{k+1,t,w}),\alpha_{t,g}(\sum\limits_{w\in \cW}e_{N-k,t,w}) \} - p_{t,g},
\end{align*}

\subsubsection*{For the lower generation limit  $\underline{\bm p}_t$ constraints}
\begin{align*}
    & \bm b_p - \bm A_p^\top \bm x = -\alpha_{t,g}\bm 1 \in \mathbb{R}^{|\cW|},\\
    & (\bm b_p \! - \! \bm A_p^\top \bm x)^\top \bm\xi_i \! + \! d_p \! - \! \bm a_p^\top \bm x =  -\underline{p}_{t,g} \! - \alpha_{t,g}(\sum\limits_{w\in \cW}e_{i,t,w}) + p_{t,g},\\ 
    & q_p(\bm x) + d_p - \bm a^\top \bm x =  \\
    & \! - \underline{p}_{t,g} \! + \! \min\{ - \alpha_{t,g}(\sum\limits_{w\in \cW}e_{k+1,t,w}), - \alpha_{t,g}(\sum\limits_{w\in \cW}e_{N-k,t,w}) \} \! + p_{t,g},
\end{align*}
Only the constraints above have the one-dimensional structure and can be strengthened as discussed in Section~\ref{sec:analysis_1dim}.
\subsubsection*{For the upper flow limit $\overline{\bm f}_t$ constraints}
{\small\begin{align*}
    & \bm b_p - \bm A_p^\top \bm x = \bm S^G_{\ell} \bm \alpha_t \bm 1 - \bm S^W_{\ell} \in \mathbb{R}^{|\cW|}, \\
    & (\bm b_p - \bm A_p^\top \bm x)^\top \bm\xi_i + d_p -\bm a_p^\top \bm x = \notag \\
    & \qquad \overline{f}_\ell - \bm S^G_{\ell} (\bm p - (\sum_{w\in \cW}e_{i,t,w})\bm\alpha_t) - \bm S^W_{\ell} (\bm \omega_t + \bm e_{i,t}) + \bm S^D_{\ell} \bm d_t.
\end{align*}}
    
\subsubsection*{For the lower flow limit $\underline{\bm f}_t$ constraints}
{\small\begin{align*}
    & \bm b_p - \bm A_p^\top \bm x = -\bm S^G_{\ell} \bm \alpha_t \bm 1 + \bm S^W_{\ell} \in \mathbb{R}^{|\cW|}, \\
    & (\bm b_p - \bm A_p^\top \bm x)^\top \bm\xi_i + d_p -\bm a_p^\top \bm x = \notag \\
    & \qquad -\underline{f}_\ell + \bm S^G_{\ell} (\bm p - (\sum_{w\in \cW}e_{i,t,w})\bm\alpha_t)  + \bm S^W_{\ell} (\bm \omega_t + \bm e_{i,t}) - \bm S^D_{\ell} \bm d_t, 
\end{align*}}where $\bm S_\ell^G$ represents the row of the PTDF matrix for line $\ell$.

\section{Proof}\label{appen:proof}

We start with introducing the following important lemma:
\begin{lemma}\label{lemma:cvar}
    $\mathcal{X}_C$ can be equivalently defined by the emperical CVaR operator $\mathbb{P}_N\mhyphen\mathrm{CVaR_{1-\epsilon}} (\cdot)$ as follows:
    \begin{align}\label{eq:ori_eq}
        \Big\{ \bm x \in \mathcal{X}: \frac{\Theta(\bm x)}{\epsilon} + \mathbb{P}_N\mhyphen\mathrm{CVaR_{1-\epsilon}}(-d(\bm x, \bm\xi)) \leq 0 \Big\},
    \end{align}
    where we have defined
    {\small
    \begin{subequations}
        \begin{align}
            & \mathbb{P}_N\mhyphen\mathrm{CVaR_{1-\epsilon}}(-d(\bm x, \bm\xi)) \coloneqq \nonumber \\
            & \quad \min_{s'} \Big\{ s' + \frac{1}{\epsilon N} \sum\limits_{i \in [N]} \max \{0, -d(\bm x, \bm\xi_i) - s'\} \Big\}, \label{eq:cvar_def} \\
            & d(\bm x, \bm\xi_i) \coloneqq \min_{p\in[P]} \Big\{ d_p(\bm x, \bm\xi_i) \Big\}, \\
            & d_p(\bm x, \bm\xi_i) \coloneqq (\bm b_p - \bm A_p^\top \bm x)^\top \bm \xi_i - \bm a_p^\top \bm x + d_p, \\
            & \Theta(\bm x) \coloneqq \theta \max\limits_{p\in[P]} \| \bm b_p - \bm A_p^\top \bm x \|_*.
        \end{align}
    \end{subequations}
    }
\end{lemma}
\begin{IEEEproof}
   \eqref{eq:ori_eq} $\Longrightarrow \mathcal{X}_C$: Similar to Lemma 2 of our previous work \cite{SFLA}, we can assume that $s'$ is selected such that the term $s' + \frac{1}{\epsilon N} \sum_{i \in [N]} \max \{0, -d(\bm x, \bm\xi_i) - s'\}$ inside the CVaR operator \eqref{eq:cvar_def} is minimised. Then, take $r_i^* = \max \{0, -d(\bm x, \bm\xi_i) - s'\}$ (implying $r_i^* \geq 0$). The inequality in~\eqref{eq:ori_eq} indicates the following:
    {\small\begin{equation}\label{eq:right_direction}
         \dfrac{\Theta(\bm x)}{\epsilon} + s' + \frac{1}{\epsilon N} \sum_{i \in [N]} r_i^* \leq 0 
         \Leftrightarrow 
         \epsilon N s' + \sum_{i \in [N]} r_i^* \leq -\Theta(\bm x) N. 
    \end{equation}  }  
    Taking $s^* = -s'$, now we will prove that the tuple $(\bm x, s^*, \bm r^*)$ satisfy all the constraints in $\mathcal{X}_C$, namely \eqref{eq:cvar-lhs-eq-s-r}--\eqref{eq:cvar-lhs-eq-main}, for any $\bm x$ belonging to the set~\eqref{eq:ori_eq}. Since $r_i^* \geq 0$ and $-\Theta(\bm x) N \leq 0$ (because of the non-negativity of dual norm),~\eqref{eq:right_direction} implies:
    \begin{align*}
        \epsilon N s^* - \sum_{i \in [N]} r_i^* \geq \Theta(\bm x) N = \theta N \max_{p\in[P]} \| \bm b_p - \bm A_p^\top \bm x \|_*,
    \end{align*}
    which means~\eqref{eq:cvar-lhs-eq-s-r} and~\eqref{eq:cvar-lhs-eq-epsNs} are satisfied. It remains to prove the satisfaction of~\eqref{eq:cvar-lhs-eq-main}. Based on the setting of $r_i^*$, we have:
    \begin{equation}
    \begin{split}
        s^* - r_i^* & = s^* - \max \Big\{0, -d(\bm x, \bm\xi_i) + s^* \Big\}.
    \end{split}
    \label{eq:s-ri}
    \end{equation}
    Notice that constraint~\eqref{eq:cvar-lhs-eq-main} is equivalent to:
    \begin{equation}
       d(\bm x, \bm\xi_i) \geq s-r_i \quad \forall i\in [N]. \label{eq:ori_main_eq}
    \end{equation}
    When $d(\bm x, \bm\xi_i) \geq s^*$,~\eqref{eq:cvar-lhs-eq-main} must be satisfied for $(\bm x, s^*, \bm r^*)$ due to $r_i^* \geq 0$. When $d(\bm x, \bm\xi_i) < s^*$, we have $s^*-r_i^* =d(\bm x, \bm\xi_i)$ based on~\eqref{eq:s-ri}. Therefore,~\eqref{eq:ori_main_eq} and equivalently~\eqref{eq:cvar-lhs-eq-main} are satisfied.

   \eqref{eq:ori_eq} $\Longleftarrow \mathcal{X}_C$: For any $\bm x$ belonging to the set $\mathcal{X}_C$, by definition there exist $(\bm r, s)$ such that constraints in the set $\mathcal{X}_C$, namely~\eqref{eq:cvar-lhs-eq-s-r}--\eqref{eq:cvar-lhs-eq-main}, are satisfied. Among these constraints,~\eqref{eq:cvar-lhs-eq-main} (equivalently~\eqref{eq:ori_main_eq}) informs that: 
    \begin{align*}
        r_i \geq s - d(\bm x, \bm\xi_i).
    \end{align*}
    Combining with $r_i \geq 0$ in~\eqref{eq:cvar-lhs-eq-s-r}, we further have:
    \begin{align*}
        r_i \geq \max \Big\{0, s - d(\bm x, \bm\xi_i) \Big\}.
    \end{align*}
    Then we have: 
    \begin{align*}
        \epsilon N s - \sum\nolimits_{i \in [N]} r_i \leq \epsilon N s  - \sum\nolimits_{i \in [N]} \max \Big\{0, s - d(\bm x, \bm\xi_i) \Big\}.
    \end{align*}
    Based on constraint \eqref{eq:cvar-lhs-eq-epsNs}, namely $\epsilon N s - \sum\nolimits_{i\in [N]}r_i \geq \Theta(\bm x) N$, we have:
    \begin{align*}
        \epsilon N s  - \sum\nolimits_{i \in [N]} \max \Big\{0, s - d(\bm x, \bm\xi_i) \Big\} \geq \Theta(\bm x) N,
    \end{align*}
    which means:
    \begin{align*}
        \frac{\Theta(\bm x)}{\epsilon} - s + \frac{1}{\epsilon N}\sum\nolimits_{i \in [N]} \max \Big\{0, s - d(\bm x, \bm\xi_i) \Big\} \leq 0.
    \end{align*}
    Now take $s' = -s$, we further have:
    \begin{align*}
        \frac{\Theta(\bm x)}{\epsilon} + s' + \frac{1}{\epsilon N}\sum\nolimits_{i \in [N]} \max \Big\{0, -s' - d(\bm x, \bm\xi_i) \Big\} \leq 0
    \end{align*}
    which leads to:
    \begin{align*}
        \frac{\Theta(\bm x)}{\epsilon} + \min_{s'}\Big\{ s' + \frac{1}{\epsilon N}\sum\nolimits_{i \in [N]} \max\{0, -s' - d(\bm x, \bm\xi_i) \Big\} \leq 0.
    \end{align*}
    Therefore, the constraint in~\eqref{eq:ori_eq} is satisfied.
\end{IEEEproof}

\vspace{-3 mm}
\subsection*{Proof of Theorem \ref{theorem:cvar_subeq_fica_exist}}
\vspace{-1 mm}

\begin{IEEEproof}
    Notice that the CVaR operator in \eqref{eq:cvar_def} can be equivalently expressed as the following form:  
    {\small
    \begin{subequations}\label{eq:CVaR_optional}
    \begin{align}
    & \mathbb{P}_N\mhyphen\mathrm{CVaR_{1-\epsilon}}(-d(\bm x, \bm\xi)) \\
    = & \min_{s'} \Big\{ s' + \frac{1}{\epsilon N} \sum\nolimits_{i \in [N]} \max \{0, -d(\bm x, \bm\xi_i) - s'\} \Big\} \\
    = & \min_{s'\in \mathbb{R}, \bm r \geq 0} \; s' + \frac{1}{\epsilon N} \sum\nolimits_{i \in [N]} r_{i} \label{eq:CVaR_optional_obj}\\
        &\qquad \mbox{s.t.} \quad  r_{i} \geq -d(\bm x, \bm\xi_i) - s',\quad i\in [N]. \notag
    \end{align}
    \end{subequations}
    }By strong duality of linear programming, we have:
    \begin{equation}\label{eq:CVaR_optional_obj_dual}
        \begin{split}
            \eqref{eq:CVaR_optional_obj} = & \max_{\bm 0\leq \bm y \leq \bm 1} \; -\frac{1}{\epsilon N} \sum\nolimits_{i\in [N]} d(\bm x, \bm\xi_i) y_{i}\\
            & \qquad \mbox{s.t.}\quad \frac{1}{N} \sum\nolimits_{i\in [N]} y_{i} = \epsilon. 
        \end{split}
    \end{equation}  
    Without loss of generality, we assume the following ordering $d(\bm x, \bm\xi_1) \leq \cdots \leq d(\bm x, \bm\xi_N)$.
    Then a primal-dual optimal pair for the problems \eqref{eq:CVaR_optional_obj} and \eqref{eq:CVaR_optional_obj_dual} is given by:
    {\small
    \begin{subequations}\label{eq:s_r_y_pri_dual}
        \begin{align}
           & s'^* = -d(\bm x, \bm\xi_{k+1}), \\
           & r_{i}^* = \max\{ 0, -d(\bm x, \bm\xi_i) - s'^* \}, \label{eq:r_pri_dual} \\
           & y_{i}^* = 
            \begin{cases} 
              1, & i = 1, \cdots, k \\
              \epsilon N - k, & i=k+1 \\
              0, & i > k+1.
           \end{cases} 
        \end{align}
    \end{subequations}       
    }The validity of the above primal-dual optimal pair holds because of the satisfaction of primal and dual constraints and the equality of the objectives of the primal and dual problems:
    \begin{align*}\label{eq:CVaR_primal_dual}
        -\frac{1}{\epsilon N} \sum\nolimits_{i\in [N]} d(\bm x, \bm\xi_i) y_{i}^*
        = s'^* + \frac{1}{\epsilon N} \sum\nolimits_{i \in [N]} r_{i}^*,
    \end{align*}
    which holds because of: 
    {\small \begin{align*}
        & -\frac{1}{\epsilon N} \sum\nolimits_{i\in [N]} d(\bm x, \bm\xi_i) y_{i}^* \\
        = &
        -\frac{1}{\epsilon N} \Big( \sum\nolimits_{i\in [k]} d(\bm x, \bm\xi_i) + d(\bm x, \bm\xi_{k+1}) (\epsilon N - k) \Big) \\
        = & -d(\bm x, \bm\xi_{k+1}) - \frac{1}{\epsilon N} \Big( \sum\nolimits_{i\in [k]} d(\bm x, \bm\xi_i) - k\cdot d(\bm x, \bm\xi_{k+1}) \Big)
    \end{align*}}and
    {\small \begin{align*}
        & s'^* + \frac{1}{\epsilon N} \sum\nolimits_{i \in [N]} r_{i}^* \\
        = & -d(\bm x, \bm\xi_{k+1}) + \frac{1}{\epsilon N} \sum_{i \in [N]} \max\{ 0, -d(\bm x, \bm\xi_i) + d(\bm x, \bm\xi_{k+1}) \} \\
        = & -d(\bm x, \bm\xi_{k+1}) + \frac{1}{\epsilon N} \sum\nolimits_{i \in [k]} \Big( -d(\bm x, \bm\xi_i) + d(\bm x, \bm\xi_{k+1}) \Big) \\
        = & -d(\bm x, \bm\xi_{k+1}) - \frac{1}{\epsilon N} \Big( \sum\nolimits_{i\in [k]} d(\bm x, \bm\xi_i) - k\cdot d(\bm x, \bm\xi_{k+1}) \Big),
    \end{align*}}where the second equality holds because of $d(\bm x, \bm\xi_{k+1})-d(\bm x, \bm\xi_i)\leq 0$ for $i\geq k+1$.
    
    Given the primal-dual analysis of the CVaR operator \eqref{eq:cvar_def} and Lemma \ref{lemma:cvar}, for any $\bm x \in \mathcal{X}_C$, we have:
    {\small
    \begin{subequations}       
        \begin{align}
            & \frac{\Theta(\bm x)}{\epsilon} + \mathbb{P}_N\mhyphen\mathrm{CVaR_{1-\epsilon}}(-d(\bm x, \bm\xi)) \leq 0 \\
            \Leftrightarrow \quad
            & \frac{\Theta(\bm x)}{\epsilon} + s'^* + \frac{1}{\epsilon N} \sum\nolimits_{i\in [N]} r_i^* \leq 0 \\
            \Leftrightarrow \quad 
            & \epsilon N s'^* + \sum\nolimits_{i\in[N]} r_i^* \leq -\Theta(\bm x) N. \label{eq:derive}
        \end{align}
    \end{subequations}   
    }Eq.~\eqref{eq:derive} informs us that $s'^* \leq 0$, due to the non-negativity of $\Theta(\bm x)$ and $r_i$ (based on \eqref{eq:r_pri_dual}).
    Now let us define $s^* \coloneqq -s'^* \geq 0$, \eqref{eq:derive} further indicates 
    {\small
    \begin{equation*}
        \epsilon N s^* - \sum\nolimits_{i \in [N]} r_i^* \geq \Theta(\bm x) N = \theta N \max\limits_{p\in[P]} \| \bm b_p - \bm A_p^\top \bm x \|_*,
    \end{equation*}
    }As $r_i^* \geq 0$, we have proved that constraints~\eqref{eq:FICA_s_r} and~\eqref{eq:FICA_epsilonNs} are satisfied. It remains to prove the satisfaction of~\eqref{eq:FICA_main}--\eqref{eq:FICA_s_str}. 
    For \eqref{eq:FICA_main} and \eqref{eq:FICA_main2}, we can even prove their stronger version, namely \eqref{eq:cvar-lhs-eq-main}. Based on the setting of $r_i^*$ in~\eqref{eq:r_pri_dual} and the definition $s^* \coloneqq -s'^* \geq 0$, we have: 
    \begin{equation}
        \label{eq:s-r}
        \begin{split}
            s^* - r_i^* & = s^* - \max \{0, -d(\bm x, \bm\xi_i) + s^*\}.
        \end{split}
    \end{equation}
    Recall that constraint~\eqref{eq:cvar-lhs-eq-main} is equivalent to the following (see~\eqref{eq:ori_main_eq}):
    \vspace{-3mm}
    \begin{align*}
       d(\bm x, \bm\xi_i) \geq s-r_i \quad i\in [N].
    \end{align*}
    When $d(\bm x, \bm\xi_i) \geq s^*$,~\eqref{eq:cvar-lhs-eq-main} must be satisfied due to $r_i^* \geq 0$. When $d(\bm x, \bm\xi_i) < s^*$, the inequality \eqref{eq:s-r} informs that $s^*-r_i^* = d(\bm x, \bm\xi_i)$, so~\eqref{eq:cvar-lhs-eq-main} is also satisfied; thus \eqref{eq:FICA_main} and \eqref{eq:FICA_main2} are satisfied. Therefore, we have proved that the pair $(\bm x, s^*, \bm r^*)$ satisfies constraints \eqref{eq:FICA_s_r}--\eqref{eq:FICA_main2}, with $s^* = -s'^* = d(\bm x, \bm\xi_{k+1})$ and $r_i^* = \max\{ 0, -d(\bm x, \bm\xi_i) + s^* \}$. Finally, based on our assumed ordering, $s^* = d(\bm x, \bm\xi_{k+1})$ is the $(k+1)$-th smallest value amongst $\{d(\bm x, \bm\xi_i)\}_{i\in [N]}$. Then for any $p\in [P]^*$, when $i\geq k+1$, we have:
    \begin{align*}
        s^*  = d(\bm x, \bm\xi_{k+1}) \leq d(\bm x, \bm\xi_i) =  \min_{p\in [P]} d_p(\bm x, \bm\xi_i)
        \leq  d_p(\bm x, \bm\xi_i).
    \end{align*}
    In other words, for any $p \in [P]^*$, there will be at least $N-k$ elements in the set $\{ d_p(\bm x, \bm\xi_i) \}_{i\in [N]}$ not less than $s^*$. Since $q_p(\bm x)$ is defined as the $(k+1)$-th smallest value of the set $\{ (\bm b_p - \bm A_p^\top \bm x)^\top \bm \xi_i \}_{i\in[N]}$, the term $q_p(\bm x) + d_p - \bm a_p^\top \bm x$ is the $(k+1)$-th smallest value of $\{ d_p(\bm x, \bm\xi_i) \}_{i\in[N]}$, and thus we must have:
    \begin{align*}
        s^* \leq q_p(\bm x) + d_p - \bm a_p^\top \bm x,\ \forall p\in[P]^*,
    \end{align*}
    which indicates the satisfaction of~\eqref{eq:FICA_s_str}.
\end{IEEEproof}

\subsection*{Proof of Theorem \ref{theorem:fica_subeq_cvar_all}}
\begin{IEEEproof}
    First, it is apparent that constraints \eqref{eq:cvar-lhs-eq-s-r} and \eqref{eq:cvar-lhs-eq-epsNs} are always satisfied for any $(\bm x, \bm r, s)$ that satisfies the FICA constraints in \eqref{eq:FICA}. It remains to prove the satisfaction of \eqref{eq:cvar-lhs-eq-main}.
    Based on the definition of $[N]_p$, we have:
    \begin{equation*}
        d_p(\bm x, \bm\xi_i) \geq q_p(\bm x) + d_p - \bm a_p^\top \bm x, \quad \forall i \in [N] \setminus [N]_p,\ p\in [P].
    \end{equation*}
    Therefore, constraint~\eqref{eq:FICA_s_str} implies that for any $i \in [N] \setminus [N]_p$ with $p\in [P]^*$, we have $d_p(\bm x, \bm\xi_i) \geq q_p(\bm x) + d_p - \bm a_p^\top \bm x \geq s$.
    As $r_i \geq 0$, we further have $d_p(\bm x, \bm\xi_i) \geq s - r_i$,
    which, combined with~\eqref{eq:FICA_main} and \eqref{eq:FICA_main2}, implies the satisfaction of~\eqref{eq:cvar-lhs-eq-main}.
\end{IEEEproof}

}

\bibliographystyle{IEEEtran}
\bibliography{references.bib}

\vfill
\end{document}